\def\mytitle{Lower Bounds for Higher-Order Convex Optimization}

\def\myauthors{Naman Agarwal\footnote{namana@cs.princeton.edu,
Computer Science, Princeton University}\hspace{2cm} Elad Hazan \footnote{ehazan@cs.princeton.edu,
Computer Science, Princeton University}}
\documentclass[11pt]{article}
\usepackage[fullpage,nousetoc,hylinks]{paper}
\usepackage{algorithm}
\usepackage{algorithmic}
\theoremstyle{plain}
%



\DeclareMathOperator*{\argmin}{argmin}
\DeclareMathOperator*{\argmax}{argmax}

\def\L{{\mathcal L}}

\def\reals{{\mathbb R}}


\newcommand{\ignore}[1]{}

\def\reals{{\mathbb R}}

\def\bold0{\mathbf{0}}





\newcommand\pr{\mbox{\bf Pr}}

\def\epsilon{\varepsilon}

\newcommand{\defeq}{\triangleq}

\newcommand{\innerprod}[2]{\langle #1, #2 \rangle}


%

\newtheorem{theorem}{Theorem}[section]

\newtheorem{claim}[theorem]{Claim}

\newtheorem{lemma}[theorem]{Lemma}
\newtheorem{corollary}[theorem]{Corollary}

\newtheorem{observation}[theorem]{Observation}

\newtheorem{definition}[theorem]{Definition}

\newcommand{\newreptheorem}[2]{%
\newenvironment{rep#1}[1]{%
 \def\rep@title{#2 \ref{##1}}%
 \begin{rep@theorem}}%
 {\end{rep@theorem}}}

\newreptheorem{theorem}{Theorem}
\newreptheorem{lemma}{Lemma}
\newreptheorem{proposition}{Proposition}
\newreptheorem{claim}{Claim}
\newreptheorem{corollary}{Corollary}
\newreptheorem{mainlemma}{Main Lemma}

\theoremstyle{remark}


\newcommand{\namedref}[2]{\mbox{\hyperref[#2]{#1~\ref*{#2}}}}

\newcommand{\sectionref}[1]{\namedref{Section}{#1}}
\newcommand{\appendixref}[1]{\namedref{Appendix}{#1}}

\newcommand{\figurerefb}[2]{\mbox{\hyperref[#1]{Figure~\ref*{#1}#2}}}

\newcommand{\equationref}[1]{\mbox{\hyperref[#1]{(\ref*{#1})}}}
\renewcommand{\eqref}{\equationref}


\numberwithin{equation}{section}


\def\smooth#1{#1_{\delta,\Gamma}}

\def\hardf{f^{\dagger}}
\def\tf{\tilde{f}}
\def\smallq{\frac{1}{20T^{1.5}}}

\BEGINDOC
\begin{titlepage}
\makeheader
\begin{abstract}
State-of-the-art methods in convex and non-convex optimization employ higher-order derivative information, either implicitly or explicitly. We explore the limitations of higher-order optimization and prove that even for convex optimization, a polynomial dependence on the approximation guarantee and higher-order smoothness parameters is necessary. As a special case, we show Nesterov's accelerated cubic regularization method to be nearly tight. 
\end{abstract}
\end{titlepage}
\pagebreak

\section{Introduction}

State-of-the-art optimization for machine learning has shifted from gradient based methods, namely stochastic gradient descent and its derivatives  \cite{adagrad,svrg}, to methods based on higher moments.  Notably, the fastest theoretical running times for both convex \cite{LiSSA2016,xu2016sub,bollapragada2016exact} and non-convex \cite{Lissa2,CarmonAGD} optimization are attained by algorithms that either explicitly or implicitly exploit second order information and third order smoothness. 

Of particular interest is Newton's method, due to recent efficient implementations that run in near-linear time in the input representation. The hope was that Newton's method, and perhaps higher order methods, can achieve iteration complexity that is independent of the condition number of the problem as well as of the dimensionality, both of which are extremely high for many large-scale applications. 

In this paper we explore the limitations of higher-order iterative optimization, and show that unfortunately, these hopes cannot be attained without stronger assumptions on the underlying optimization problem.  To the best of our knowledge, our results are the first  lower bound for $k^{th}$ order optimization for $k \geq 2$ that includes higher-order smoothness. \footnote{After the writing of the first manuscript we were made aware of the work by Arjevani et al. \cite{shamir2017oracle} which provides lower bounds for these settings as well.}


We consider the problem of $k^{th}$-order optimization. We model a $k^{th}$-order algorithm as follows. Given a $k$-times differentiable function $f: \reals^d \rightarrow R$, at every iteration $i$ the algorithms outputs a point $x_i$ and receives as input the tuple $[ f(x_i), \nabla f(x_i), \nabla^2 f(x_i) \ldots \nabla^k f(x_i)]$, i.e. the value of the function and its $k$ derivatives at $x$ \footnote{An iteration is equivalent to an oracle call to $k^{th}$-order derivatives in this model}. The goal of the algorithm is to output a point $x_T$ such that 
\[ f(x_T) - \min_{x \in \reals^d} f(x) \leq \epsilon \]
For the $k^{th}$-order derivatives to be informative, one needs to bound their rate of change, or Lipschitz constant. This is called $k^{th}$-order smoothness, and we denote it by $L_{k+1}$. In particular we assume that
$$ \| \nabla^k f(x) - \nabla^k f(y) \| \leq L_{k+1} \| x-y \| $$
where $\|\nabla^k f\|$ is defined as the induced operator norm with respect to the Euclidean norm. Our main theorem shows the following limitation of $k^{th}$-order iterative optimization algorithms: 

\begin{theorem}
\label{thm:mainthminit}
For every number $\L_{k+1}$ and $k^{th}$-order algorithm ALG (deterministic or randomized), there exists an $\epsilon_0(\L_{k+1})$ such that for all $\epsilon \leq \epsilon_0(\L_k)$, there exists a $k$-differentiable convex function $f \in B_d \rightarrow \reals$ with $k^{th}$-order smoothness coefficient $\L_{k+1}$ such that ALG cannot output a point $x_T$ such that
\[f(x_T) \leq \min_{x \in \reals^d} f(x) + \epsilon\]
 in number of iterations $T$ fewer than
\[ c_k\left(\frac{\L_{k+1}}{\epsilon}\right)^{2/(5k + 1)}\]
where $c_k$ is a constant depending on $k$ and $B_d$ is defined to be the unit ball in $d$ dimensions.
\end{theorem}

Note that although the bound is stated for constrained optimization over the unit ball, it can be extended to an unconstrained setting via the addition of an appropriate scaled multiple of $\|x\|^2$. We leave this adaptation for a full version of this paper. Further as is common with lower bounds the underlying dimension $d$ is assumed to be large enough and differs for the determinisitc vs randomized version. Theorems \ref{thm:mainthm} and \ref{thm:mainthmrandomized} make the dependence precise.

\paragraph{Comparison to existing bounds. } 
For the case of $k=2$, the most efficient methods known are the cubic regularization technique proposed by Nesterov \cite{nesterov2008cubic} and an accelerated hybrid proximal extragradient method proposed by Monteiro and Svaiter \cite{monteiro2013accelerated}. The best known upper bound in this setting is $O\left(\frac{L_3}{\epsilon}\right)^{2/7}$\cite{monteiro2013accelerated}. We show a lower bound of $\Omega\left(\left( \frac{L_3}{\epsilon}\right)^{2/11}\right)$ demonstrating that the upper bound is nearly tight. 

For the case of case of $k > 2$, Baes \cite{baeshigherorder} proves an upper bound of $O\left( \left(\frac{L_{k+1}}{\epsilon_k} \right)^{\frac{1}{k+1}}\right)$. In comparison Theorem \ref{thm:mainthminit} proves a lower bound of $\Omega\left(\left(\frac{L_{k+1}}{\epsilon}\right)^{2/(5k + 1)}\right)$. 

\subsection{Related work}
The literature on convex optimization is too vast to survey; the reader is referred to \cite{boyd,NesterovBook}. 

Lower bounds for convex optimization were studied extensively in the seminal work of \cite{Nemirovsky1978}. In particular, tight first-order optimization lower bounds were established assuming first-order smoothness.(Also see \cite{NesterovBook} for a concise presentation of the lower bound). In a recent work \cite{arjevanisecondorder} presented a lower bound when given access to second-order derivatives. However a key component (as remarked by the authors themselves) missing from the bound established by \cite{arjevanisecondorder} was that the constructed function was not third-order smooth. Indeed the lower bound established by \cite{arjevanisecondorder} can  be overcome when the function is third-order smooth (ref. \cite{nesterov2008cubic}). The upper bounds for higher-order oracles (assuming appropriate smoothness) was established by \cite{baeshigherorder}. 

Higher order smoothness has been leveraged recently in the context of non-convex optimization \cite{Lissa2,CarmonAGD,allen2017natasha}. In a surprising new discovery, \cite{carmon2017convex} show that assuming higher-order smoothness the bounds for first-order optimization can be improved without having explicit access to higher-order oracles. This is a property observed in our lower bound too. Indeed as shown in the proof the higher order derivatives at the points queried by the algorithm are always 0. For further details regarding first-order lower bounds for various different settings we refer the reader to \cite{agarwal2009information,woodworth2016tight,arjevanisecondorder,arjevani2015lower} and the references therein. 

In parallel and independently, Arjevani et al. \cite{shamir2017oracle} also obtain lower bounds for deterministic higher-order optimization.  In comparison, their lower bound is stronger in terms of the exponent than the ones proved in this paper, and matches the upper bound for $k=2$.  However, our construction and proof are simple (based on the well known technique of smoothing) and our bounds hold for randomized algorithms as well, as opposed to the their deterministic lower bounds.   


\subsection{Overview of Techniques}

Our lower bound is inspired by the lower bound presented in \cite{ClarksonHW2012}. In particular we construct the function as a piecewise linear convex function defined by $f(x) = \max_i \{ a_i^Tx \}$ with carefully constructed vectors $a_i$ and restricting the domain to be the unit ball. The key idea here is that querying a point reveals information about at most one hyperplane. The optimal point however can be shown to require information about all the hyperplanes.

Unfortunately the above function is not differentiable. We now smooth the function by the ball smoothing operator (defined in Definition \ref{def:smoothingoperator}) which averages the function in a small Euclidean ball around a point. We show (c.f. Corollary \ref{cor:maincor}) that iterative application of the smoothing operator ensures $k$-differentiability as well as boundedness of the $k^{th}$-order smoothness.

Two key issues arise due to above smoothing. Firstly although the smoothing operator leaves the function unchanged around regions far away from the intersection of the hyperplanes, it is not the case for points lying near the intersection. Indeed querying a point near the intersection of the hyperplanes can potentially lead to leak of information about multiple hyperplanes at once. To avoid this, we carefully shift the linear hyperplanes making them affine and then arguing that this shifting indeed forces sufficient gap between the points queried by the algorithm and the intersections leaving sufficient room for smoothing. 

Secondly such a smoothing is well known to introduce a dependence on the dimension $d$ in the smoothness coefficients. Our key insight here is that for the class of functions being considered for the lower bound (c.f. Definition \ref{defn:gammainvariance}) smoothing can be achieved without a dependence on the dimension(c.f. Theorem \ref{thm:smoothingmain}). This is essential to achieving dimension free lower bounds and we believe this characterization can be of intrinsic interest.

\subsection{Organization of the paper}

We begin by providing requisite notations and definitions for the smoothing operator and proving the relevant lemmas regarding smoothing in Section \ref{sec:prelimssmooth}. In Section \ref{sec:construction} we provide the construction of our hard function. In Section \ref{sec:deterministic} we state and prove our main theorem (Theorem \ref{thm:mainthm}) showing the lower bound against determinsitic algorithms. We also prove Theorem \ref{thm:mainthminit} based on Theorem \ref{thm:mainthm} in this Section. In Section \ref{sec:randomized} we state and prove the Theorem \ref{thm:mainthmrandomized} showing the lower bound against randomized algorithms.

\section{Preliminaries}
\label{sec:prelimssmooth}
\subsection{Notation}

We use $B_d$ to refer to the $d$ dimensional $\ell_2$ unit ball. We suppress the $d$ from the notation when it is clear from the context. Let $\Gamma$ be an $r$-dimensional linear subspace of $\reals^d$. We denote by $M_{\Gamma}$, an $r \times d$ matrix which containes an orthonormal basis of $\Gamma$ as rows. Let $\Gamma^{\perp}$ denote the orthogonal complement of $\Gamma$. Given a vector $v$ and a subspace $\Gamma$ let $v \perp \Gamma$ denote the perpendicular component of $v$ w.r.t $\Gamma$. We now define the notion of a $\Gamma$ invariant function. 

\begin{definition}[$\Gamma$-invariance]
\label{defn:gammainvariance}
Let $\Gamma$ be an $r$ dimensional linear subspace of $\reals^d$. A function $f:\reals^d \rightarrow \reals$ is said to be $\Gamma$-invariant if for all $x \in \reals^d$ and $y$ belonging to the subspace $\Gamma^{\perp}$, i.e. $M_{\Gamma} y = 0$ we have that
    \[ f(x) = f(x + y)\]
   Equivalently there exists a function $g:\reals^r \rightarrow \reals$ such that for all $x$, $f(x) = g(M_{\Gamma} x)$.
\end{definition}
\noindent A function $f: \reals^d \rightarrow \reals$ is defined to be $c$-Lipschitz with respect to a norm $\|\cdot\|$ if it satisfies
\[ f(x) - f(y) \leq c\|x-y\| \] 
 Lipschitzness for the rest of the paper will be measured in the $\ell_2$ norm. 
\subsection{Smoothing}
In this section we define the smoothing operator and the requisite properties.
\begin{definition}[Smoothing operator]
\label{def:smoothingoperator}
Given a $r$-dimensional subspace $\Gamma \in \reals^d$ and a parameter $\delta > 0$ define the operator $S_{\delta,\Gamma}:(\reals^d \rightarrow \reals) \rightarrow (\reals^d \rightarrow \reals)$ (refered henceforth as the smoothing operator) as 
  \[ S_{\delta, \Gamma}f(x) \defeq E_{v \in \Gamma, \|v\| \leq 1}[ f(x + \delta v)]\]
  As a shorthand we define $f_{\delta, \Gamma} \defeq S_{\delta, \Gamma}f$. Further for any $t \in \mathbb{N}$ define $S_{\delta,\Gamma}^t f \defeq S_{\delta,\Gamma}(... S_{\delta,\Gamma}( f )) $ i.e. the smoothing operator applied on $f$ iteratively $t$ times.
\end{definition}
When $\Gamma = \reals^d$ we suppress the notation from $f_{\delta, \Gamma}$ to $f_{\delta}$. Following is the main lemma we prove regarding the smoothing operator. 
\begin{lemma}
  \label{thm:smoothingmain}
  Let $\Gamma $ be a $r$ dimensional linear subspace of $\reals^d$ and $f : \reals^d \rightarrow \reals$ be $\Gamma$-invariant and $G$-Lipschitz. Let $\smooth{f} \defeq S_{\delta, \Gamma} f$ be the smoothing of $f$. Then we have the following properties.
  \begin{enumerate}
    \item $\smooth{f}$ is differentiable and also $G$-lipschitz and $\Gamma$-invariant.
    \item $\nabla \smooth{f}$ is $\frac{rG}{\delta}$-Lipschitz.
    \item $\forall \;x: |f_{\delta,\Gamma}(x) - f(x)| \leq \delta G$
  \end{enumerate}
  \end{lemma}
\begin{proof}
  As stated before $f$ being $\Gamma$-invariant implies that there exists a function $g$ such that $f(x) = g(\Gamma x)$. Therefore we have that
  \[ f_{\delta, \Gamma}(x) = E_{v \in \Gamma, \|v\| \leq 1}[ f(x + \delta v)] = E_{v \in \Gamma, \|v\| \leq 1}[ g(M_{\Gamma} x + \delta M_{\Gamma} v)] = g_{\delta}(M_{\Gamma} x)\]
  where $g_{\delta}(x) \defeq S_{\delta, \reals^r} g(x)$. The representation of $\smooth{f}$ as $g_{\delta}(M_{\Gamma} x)$ implies that $f_{\delta, \Gamma}$ is $\Gamma$ invariant. Further the above equality implies that $\nabla \smooth{f}(x) = M_{\Gamma}^T \nabla g_{\delta}(M_{\Gamma} x)$. A standard argument using Stokes' theorem shows that $g_{\delta}$ is differentiable even when $g$ is not \footnote{We need $g$ to be not differentiable in a measure 0 set which is always the case with our constructions} and that $\nabla g_{\delta}(y) = \frac{r}{\delta}E_{v \sim S_r} \left[ g(y + \delta v) v\right]$(Lemma 1 \cite{flaxman2005online}),
  where $S_r$ is the $r$-dimensional sphere, i.e. $S_r = \{x \in \reals^r | \|x\| = 1\}$
  \begin{align*}
    \| \nabla g_{\delta}(x) - \nabla g_{\delta}(y) \| &=  \frac{r}{\delta}\|  E_{v \sim S_r} \left[ g(x + \delta v) v\right] - E_{v \sim S_r} \left[ g(y + \delta v) v\right] \|  \\ 
    &\leq   \frac{r}{\delta} E_{v \sim S_r} \left[ | g(x + \delta v)  -  g(y + \delta v) | \| v \| \right] \leq \frac{rG}{\delta} \|x - y\| 
  \end{align*}

  The first inequality follows from Jensen's inequality and the second inequality follows from noticing that $f$ being $G$-Lipschitz implies that $g$ is $G$-Lipschitz. We now have that 
  \[ \|\nabla \smooth{f}(x) - \nabla \smooth{f}(y)\| \leq \|\nabla g_{\delta}(M_{\Gamma} x) - \nabla g_{\delta}(M_{\Gamma} y) \| \leq \frac{r G}{\delta} \|M_{\Gamma} (x - y)\| \leq \frac{r G}{\delta} \|(x - y)\|\]
  $f$ being $G$-Lipschitz immediately gives us $\forall \;x: |f_{\delta,\Gamma}(x) - f(x)| \leq \delta G$
\end{proof}

\begin{corollary}
\label{cor:maincor}
  Given a $G$-Lipschitz continuous function $f$ and an $r$-dimensional subspace $\Gamma$ such that $f$ is $\Gamma$-invariant, we have that the function $S_{\delta,\Gamma}^k f$ is $k$-times differentiable $\forall\;k$. Moreover we have that for any $x,y$
  \[ \forall i \in [k] \;\; \| \nabla^{i} S_{\delta,\Gamma}^k f(x) - \nabla^{i} S_{\delta,\Gamma}^k f(y)\| \leq \left(\frac{r}{\delta} \right)^iG\|x-y\|\] 
  \[| S_{\delta,\Gamma}^k f(x) - f(x)| \leq G\delta k \]
\end{corollary}
\begin{proof}
We will argue inductively. The base case ($k=0$) is a direct consequence of the function $f$ being $G$-Lipschitz. Suppose the theorem holds for $k-1$. To argue about $\|\nabla^{k} S_{\delta,\Gamma}^k f(x) - \nabla^{k} S_{\delta,\Gamma}^k f(y)\|$ we will consider the function $q_{i,v}(x) = \nabla^{i} S_{\delta,\Gamma}^{k-1} f(x)[v^{\otimes i}]$ for $i \in [k]$ and for a unit vector $v$. 
We will first consider the case $i < k$. Using the inductive hypothesis and the fact that smoothing and derivative commute for differentiable functions we have that 
\[S_{\delta, \Gamma}\;q_{i,v}(x) = \nabla^i S_{\delta, \Gamma}^k  f(x)[v^{\otimes i}]\]
Note that the inductive hypothesis implies that $q_{i,v}(x)$ is $\left(\frac{r}{\delta}\right)^{i}G$-Lipschitz and so is $S_{\delta, \Gamma}\;q_{i,v}(x)$ via Lemma \ref{thm:smoothingmain}. Therefore we have that 
\[ \forall i \in [k-1] \;\; \| \nabla^{i} S_{\delta,\Gamma}^k f(x) - \nabla^{i} S_{\delta,\Gamma}^k f(y)\| \leq \left(\frac{r}{\delta} \right)^iG\|x-y\|\]
We now consider the case when $i = k$. By Lemma \ref{thm:smoothingmain} we know that $S_{\delta, \Gamma}\;q_{i,v}(x) = \nabla^i S_{\delta, \Gamma}^k  f(x)[v^{\otimes i}]$ is differentiable and therefore we have that $S_{\delta, \Gamma}^k  f(x)$ is $k$ times differentiable. Further we have that 
\[\nabla S_{\delta, \Gamma} \; q_{k,v}(x) = \nabla^{k} S_{\delta,\Gamma}^k f(x)[v^{\otimes k-1}]\]
A direct application of Lemma $\ref{thm:smoothingmain}$ gives that 
\[ \|\nabla^{k} S_{\delta,\Gamma}^k f(x) - \nabla^{k} S_{\delta,\Gamma}^k f(y)\| \leq \left(\frac{r}{\delta}\right)^{k}G\|x-y\|\]
Further it is immediate to see that 
\[ \inf_{y: \|y - x\| \leq k\delta} f(y) \leq S_{\delta, \Gamma}^k f(x) \leq \sup_{y: \|y - x\| \leq k\delta} f(y)\]
which implies using the fact that $f$ is $G$ Lipschitz that 
\[| S_{\delta,\Gamma}^k f(x) - S_{\delta,\Gamma}^k f(x)| \leq G\delta k \]

\end{proof}

\section{Construction of the hard function}

\label{sec:construction}

In this section we describe the construction of our hard function $\hardf$. Our construction is inspired by the information-theoretic hard instance of zero sum games of \cite{clarkson2012sublinear}. The construction of the function will be characterized by a sequence of vectors $X^{1 \rightarrow r} = \{x_1 \ldots x_r\}$, $x_i \in B_d$ and parameters $k,\gamma, \delta, m$. We assume $d > m \geq r$. To make the dependence explicit we denote the hard function as  \[\hardf(X^{1 \rightarrow r},\gamma,k,\delta, m):B_d \rightarrow \reals\]

For brevity in the rest of the section we supress $X^{1 \rightarrow r}, \gamma, k, \delta, m$ from the notation, however all the quantities defined in the section depend on them. To define $\hardf$ we will define auxilliary vectors $\{a_1 \ldots a_r\}$ and auxilliary functions $f,\tilde{f}$. 

Given a sequence of vectors $\{x_1, x_2, \ldots x_r\}, x_i \in B_d$, let $X_i$ for $i \leq r$ be defined as the subspace spanned by the vectors $\{x_1 \ldots x_{i}\}$. Further inductively define vectors $\{a_1 \ldots a_r\}$ as follows.

\noindent If $x_i \perp X_{i-1} \neq 0$ (i.e. the perpendicular component of $x_i$ on the subspace $X_{i-1}$ is not zero), define  
\[a_i \defeq \frac{\hat{a}_i}{\| \hat{a}_i\|}\text{ where } \hat{a}_i \defeq x_i \perp X_{i-1}\]
If indeed $x_i$ belongs to the subspace $X_i$, then $a_i$ is defined to be an arbitrary unit vector in the perpendicular subspace $X_{i-1}^{\perp}$. Further define an auxilliary function 
\[f(x) \defeq \max_{i \in [r]} f_i(x) \text{ where } f_i(x) \defeq a_i^T x\]
Given the parameter $\gamma$, now define the following functions 
\[ \tilde{f}(x) \defeq \max_i \tilde{f}_i(x) \text{ where } \tilde{f}_i(x) \defeq f_i(x) + \left(1 - \frac{i}{m} \right)\gamma \defeq a_i^T x + \left(1 - \frac{i}{m} \right)\gamma\]
With these definitions in place we can now define the hard function parameterized by $k, \delta$. Let $A$ be the subspace spanned by $\{a_1 \ldots a_r\}$
\begin{equation}
  \label{eqn:hardfdef}
  \hardf(X^{1 \rightarrow r}, k, \gamma, \delta, m) \defeq S^k_{\delta, A} \;\tilde{f}(X^{1 \rightarrow r}, \gamma, m)
\end{equation}
i.e. $\hardf$ is constructed by smoothing $\tilde{f}$ $k$-times with respect to the parameters $\delta, A$. We now collect some important observations regarding the function $\hardf$. 

\begin{observation}
	$\hardf$ is convex and continuous. Moreover it is 1-Lipschitz and is invariant with the respect to the r dimensional subspace $A$.
\end{observation}

Note that $\tilde{f}$ is a $\max$ of linear functions and hence convex. Since smoothing preserves convexity we have that $\hardf$ is convex. 1-Lipschitzness follows by noting that by definition $\|a_i\| = 1$. It can be seen that $\tilde{f}$ is $A$-invariant and therefore by Theorem \ref{thm:smoothingmain} we get that $\hardf$ is $A$-invariant.

\begin{observation}
	$\hardf$ is $k$-differentiable with the Lipschitz constants $L_{i+1} \leq \left(\frac{r}{\delta}\right)^i$ for all $i \leq k$.
\end{observation}
\noindent Above is a direct consequence of Corollary \ref{cor:maincor} and the fact that $\tilde{f}$ is 1-Lipschitz and invariant with respect to the $r$-dimensional subspace $A$. 
Corollary \ref{cor:maincor} also implies that 
\begin{equation}
  \label{eqn:errorbound}
  \forall x \;\;|\hardf(x) - \tf(x)| \leq k\delta 
\end{equation}
Setting $\hat{x} \defeq -\sum \frac{a_i}{\sqrt{r}}$, we get that $f(\hat{x}) = \frac{-1}{\sqrt{r}}$. Therefore the following inequality follows from Equation \eqref{eqn:errorbound} and by noting that $\|f(x) - \tf(x)\|_{\infty} \leq \gamma$ 

\begin{equation}
\label{eqn:funcminbound}
\min_{x \in B_d} \hardf(x) \leq \hardf(\hat{x}) \leq \frac{-1}{\sqrt{r}} + \gamma + k\delta
\end{equation}
The following lemma provides a characterization of the derivatives of $\hardf$ at the points $x_i$. 

\begin{lemma}
\label{lemma:derivatives}
  Given a sequence of vectors $\{x_1 \ldots x_r\}$ and parameters $\delta, \gamma, r, m$, let $\{g_1 \ldots g_r\}$ be a sequence of functions defined as
    \[ \forall\;i\;\;g_i \defeq \hardf(X^{1 \rightarrow i}, k,\gamma, \delta, m)\] 
	If the parameters are such that 
  \[2k\delta \leq \frac{\gamma}{m}\]
  then we have that 
  \[ \forall\;i\in[r] \;\forall j \in [k]  \;\; g_i(x_i) = g_r(x_i), \; \nabla^j g_i(x_i) = \nabla^j g_r(x_i)  \]

\end{lemma}

\begin{proof}

We will first note the following about the smoothing operator $S^k_{\delta}$. At any point $x$ all the $k$ derivatives and the function value of $S^k_{\delta} f$ for any function $f$ depend only on the values of the function $f$ in a ball of radius atmost $k \delta$ around the point $x$. Consider the function $g_r$ and $g_i$ for any $i \in [r]$. Note that by definition of the functions $g_i$, for any $x$ such that \[
\argmax_{j \in [r]} \;\;a_j^Tx + \left(1 - \frac{j}{m}\right)\gamma \leq i\]
we have that $g_i(x) = g_r(x)$. Therefore to prove the lemma it is sufficient to show that
\[ \forall \;i, x \in \|x - x_i\| \leq k\delta\;\;\;\; \argmax_{j \in [r]} \;\;a_j^Tx + \left(1 - \frac{j}{m}\right)\gamma \leq i\]
Lets first note the following facts. By construction we have that $\forall j > i, a_j^Tx_i = 0$. This immediately implies that  
\begin{equation}
\max_{j > i} \;\; a_j^Tx + \left(1 - \frac{j}{m}\right)\gamma =\left(1 - \frac{i+1}{m} \right) \gamma
\end{equation}
Further using the fact that $\|a_j\| \leq 1$, $\forall j \in [r]$ we have that 
\begin{equation}
\forall x \;\;\text{s.t.}\;\;\|x - x_i\| \leq k\delta \;\;\text{we have}\;\; \max_{j > i} \;\; a_j^Tx_i + \left(1 - \frac{j}{m}\right)\gamma \leq\left(1 - \frac{i+1}{m} \right) \gamma + k\delta
\end{equation}
Further note that by construction $a_i^Tx_i \geq 0$ which implies $a_i^Tx + \left(1 - \frac{j}{m}\right)\gamma \geq \left(1 - \frac{i}{m} \right)\gamma$. Again using the fact that $\|a_j\| \leq 1$, $\forall j \in [r]$ we have that 
\begin{equation}
\label{eqn:templabel1}
\forall x \;\;\text{s.t.}\;\;\|x - x_i\| \leq k\delta \;\;\text{we have}\;\; \max_{j \leq i} \;\; a_j^Tx + \left(1 - \frac{j}{m}\right) \geq \left(1 - \frac{i}{m} \right) \gamma - k\delta
\end{equation}
The above equations in particular imply that as long as $2k\delta < \frac{\gamma}{m}$ , we have that  
\begin{equation}
\forall x \;\;\text{s.t.}\;\;\|x - x_i\| \leq k\delta \;\; \argmax_{j \in [r]} a_j^Tx + \left(1 - \frac{j}{m}\right) \leq i 
\end{equation}
which as we argued before is sufficient to prove the lemma.
\end{proof}

\section{Main Theorem and Proof}
\label{sec:deterministic}
The main theorem (Theorem \ref{thm:mainthm}) follows immediately from the following main technical lemma by setting $\epsilon = \frac{1}{2\sqrt{T}}$. \footnote{For the randomized version the statement follows the same way from Theorem \ref{thm:mainthmrandomized}}.

\begin{theorem}
\label{thm:mainthm}
  For any integer $k$, any $T > 5k$, and $d > T$ and any $k$-order deterministic algorithm working on $\reals^d$, there exists a convex function $\hardf: B_d \rightarrow \reals$ for $d > T$, such that for $T$ steps of the algorithm every point $y \in B_d$ queried by the algorithm is such that 
  \[\hardf(y) \geq \min_{x \in B_d}\hardf(x) + \frac{1}{2\sqrt{T}}.\] 
  Moreover the function is guaranteed to be $k$-differentiable with Lipschitz constants $L_{i+1}$ bounded as 
  \begin{equation}
  \label{eqn:condlk}
  	\forall \;i \leq k \;\; L_{i+1} \leq (10k)^iT^{2.5i}
  \end{equation}

\end{theorem}

We first prove Theorem \ref{thm:mainthminit} in the deterministic case using Theorem \ref{thm:mainthm}. 

\begin{proof}[Proof of Theorem \ref{thm:mainthminit} Deterministic case]
Given an algorithm ALG and numbers $\L_{k+1},k$ define $\epsilon_0(\L_{k+1},k) \defeq \L_{k+1}/(10k)^k$. For any $\epsilon \leq \epsilon_0$ pick a number $\mathcal{T}$ such that
\[\epsilon = \frac{\L_{k+1}}{(10k)^k\mathcal{T}^{(2.5k + 0.5)}}\]
Let $\hardf$ be the function constructed in Theorem \ref{thm:mainthm} for parameters $k, \mathcal{T}, ALG$ and define the hard function $h:B_d \rightarrow \reals$
\[ h(x) \defeq \frac{\L_{k+1}}{(10k)^k\mathcal{T}^{2.5k}}\hardf(x)\]
Note that by the guarantee in Equation \eqref{eqn:condlk} we get that $h(x)$ is $k^{th}$-order smooth with coefficient at most $\L_{k+1}$. Note that since this is a scaling of the original hard function $\hardf$ the lower bound applies directly and therefore ALG cannot achieve accuracy 
\[ \frac{\L_{k+1}}{(10k)^k\mathcal{T}^{2.5k}2\sqrt{T}} \defeq \epsilon \]
in less that $T = c_{k}\left(\frac{\L_{k+1}}{\epsilon}\right)^{\frac{2}{5k+1}}$ iterations which finishes the proof of the theorem.  
\end{proof}
We now provide the proof of Theorem \ref{thm:mainthm}.
\begin{proof}[Proof of Theorem \ref{thm:mainthm}]

Define the following parameters
\begin{equation}
	\gamma \defeq \frac{1}{3\sqrt{T}} \qquad\delta_T \defeq \frac{\gamma}{3kT}
\end{equation}

Consider a deterministic algorithm $Alg$. Since $Alg$ is deterministic let the first point played by the algorithm be fixed to be $x_1$. We now define a series of functions $\hardf_i$ inductively for all $i = \{1, \ldots T\}$ as follows
\begin{equation}
	\label{eqn:Xdefn}
	X^{1 \rightarrow i} \defeq \{x_1 \ldots x_i\}
\end{equation}

\begin{equation}
	\label{eqn:inducthardfdefn}
	\hardf_i \defeq \frac{\hardf(X^{1 \rightarrow i}, \gamma, k, \delta_T, T)}{\max_{x \in B_d} \hardf(X^{1 \rightarrow i}, \gamma, k, \delta_T, T)(x)}
\end{equation}

\begin{equation}
	\label{eqn:inductinpidefn}
	Inp^x_i \defeq \{ \hardf_i(x_i), \nabla \hardf_i(x_i) \ldots \nabla^k \hardf_i(x_i)\}
\end{equation}

\begin{equation}
	\label{eqn:xiindictdefn}
	x_{i+1} \defeq Alg(Inp^x_0, \ldots Inp^x_i)
\end{equation}

The above definitions \textit{simulate} the deterministic algorithm Alg with respect to changing functions $\hardf_i$. $Inp_{i}^x$ is the input the algorithm will receive if it queried point $x_i$ and the function was $\hardf_i$. $x_{i+1}$ is the next point the algorithm Alg will query on round $i+1$ given the inputs $\{Inp_1^x \ldots Inp_i^x\}$ over the previous rounds. Note that thus far these quantities are tools defined for analysis. Since Alg is deterministic these quantities are all deterministic and well defined. We will now prove that the function $\hardf_T$ defined in the series above satisfies the properties required by the Theorem \ref{thm:mainthm}. 
\\
\\
\noindent \textbf{Bounded Lipschitz Constants} Using Corollary \ref{cor:maincor} and the fact that $\hardf$ has Lipschitz constant bounded by 1 we get that the function $\hardf_T$ has higher order Lipschitz constants bounded above as
\[ \forall i \leq k \;\;\; L_{i+1} \leq \left(\frac{T}{\delta_T}\right)^i \leq \left(10kT^{2.5}\right)^i\] 
\noindent \textbf{Suboptimality}

Let $\{y_0 \ldots y_T\}$ be the points queried by the algorithm Alg when executed on $\hardf_T$. We need to show that 
\begin{equation}
	\label{eqn:suboptimalitygoal}
	\forall i \in [1 \ldots T] \qquad \hardf_T(y_i) \geq \min_{x \in B_d} \hardf_T(x) + \frac{1}{2\sqrt{T}}
\end{equation}
Equation \ref{eqn:suboptimalitygoal} follows as a direct consequence of the following two claims. 

\begin{claim}
\label{claim:consistency}
	We have that for all $i \in [1,T]$, $y_i = x_i$ where $x_i$ is defined by Equation \ref{eqn:xiindictdefn}.
\end{claim}

\begin{claim}
\label{claim:xsuboptimality}
We have that
\[  \forall i \in [1 \ldots T] \qquad \hardf_T(x_i) \geq \min_{x \in B_d} \hardf_T(x) + \frac{1}{2\sqrt{T}} \]
\end{claim} 

To remind the reader $x_i$ were variables defined by Equation \eqref{eqn:xiindictdefn} and $y_i$ are the points played by the algorithm Alg when run on $\hardf_T$. Claim \ref{claim:consistency} shows that even though $\hardf_T$ was constructed using $x_i$ the outputs produced by the algorithm does not change. 

Claim \ref{claim:consistency} and Claim \ref{claim:xsuboptimality} derive Equation \ref{eqn:suboptimalitygoal} in a straightforward manner thus finishing the proof of Theorem \ref{thm:mainthm}.

\end{proof}

We now provide the proofs of Claim \ref{claim:consistency} and Claim \ref{claim:xsuboptimality}. 

\begin{proof}[Proof of Claim \ref{claim:consistency}]

We will prove the claim inductively. The base case $x_1 = y_1$ is immediate because Alg is deterministic and therefore the first point queried by it is always the same. Further note that $y_i$ for $i \geq 2$ is defined inductively as follows. 
\begin{equation}
	Inp^y_i \defeq \{ \hardf_T(y_i), \nabla \hardf_T(y_i), \ldots \nabla^k \hardf_T(y_i)\}
\end{equation}
\begin{equation}
  \label{eqn:inductivecaseproof}
	y_{i+1} = Alg(Inp^y_1, \ldots Inp^y_T)
\end{equation}
It is now sufficient to show that 
\begin{equation}
\label{eqn:inputequalities}
  \{\forall j \leq i \;\; x_j = y_j \} \Rightarrow Inp_i^y = Inp_i^x.	
  \end{equation}
where $Inp_i^x$ is as defined in Equation \eqref{eqn:inductinpidefn}. To see this note that 
\[\{\forall j \leq i \;\; x_j = y_j \} \Rightarrow \{\forall j \leq i \;\; Inp_j^y = Inp_j^x \} \Rightarrow x_{i+1} = y_{i+1}\]
Equation \eqref{eqn:inputequalities} is a direct consequence of Lemma \ref{lemma:derivatives} by noting that $2k\delta_T \leq \gamma/T$ which is true by definition of these parameters.
\end{proof}

\begin{proof}[Proof of Claim \ref{claim:xsuboptimality}]
  Using Lemma \ref{lemma:derivatives} we have that $\hardf_i(x_i) = \hardf_T(x_i)$. Further Equation \eqref{eqn:templabel1} implies that 
  \[\hardf_i(x_i) \geq \frac{(1 - \frac{i}{T})\gamma - k\delta_T}{1 + \gamma - k\delta_T}\]Now using \eqref{eqn:funcminbound} using we get that every point in $\{x_1 \ldots x_T\}$ is such that 
\[ \hardf_T(x_i) - \min_{x \in B} \hardf_T(x) \geq  \frac{\left( \frac{1}{\sqrt{T}} - \frac{i\gamma}{T} - 2k\delta_T \right)}{1 + \gamma - k\delta_T} \geq \frac{1}{2\sqrt{T}} \] 
The above follows by the choice of parameters and $T$ being large enough. This finishes the proof of Claim \ref{claim:xsuboptimality}.
\end{proof}

\section{Lower Bounds against Randomized Algorithms}

\label{sec:randomized}

In this section we prove the version of Theorem \ref{thm:mainthm} for randomized algorithms. The key idea underlying the proof remains the same. However since we cannot \textit{simulate} the algorithm anymore we choose the vectors $\{a_i\}$ forming the subspace randomly from $\reals^d$ a large enough $d$. This ensures that no algoirthm with few queries can discover the subspace in which the function is non-invariant with reasonable probability. Naturally the dimension required for Theorem \ref{thm:mainthm} now is larger than the tight $d > T$ we achieved as in the case of deterministic algorithms. 

\begin{theorem}
\label{thm:mainthmrandomized}
  For any integer $k$, any $T > 5k$, $\delta \in [0,1]$,  and any $k$-order (potentially randomized algorithm), there exists a differentiable convex function $\hardf: \reals^d \rightarrow \reals$ for $d = \Omega(T^3\log(\delta T^2))$, such that with probability at least $1 - \delta$ (over the randomness of the algorithm) for $T$ steps of the algorithm every point queried by the algorithm is at least $\frac{1}{2\sqrt{T}}$ sub-optimal over the unit ball. Moreover the function $\hardf$ is guaranteed to be $k$-differentiable with Lipschitz constants $L_i$ bounded as 
  \[ \forall \;i \leq k \;\; L_{i+1} \leq (20kT^{2.5})^i\]
\end{theorem}
\begin{proof}
  We provide a randomized construction for the function $\hardf$.
  The construction is the same as in Section \ref{sec:construction} but we repeat it here for clarity. 
  We sample a random $T$ dimensional basis $\{a_1 \ldots a_T\}$. Let $A_i$ be the subspace spanned by $\{a_1 \ldots a_i\}$ and $A_i^{\perp}$ be the perpendicular subspace. Further define an auxilliary function 
\[f(x) \defeq \max_i f_i(x) \text{ where } f_i(x) \defeq a_i^T x\]
Given a parameter $\gamma$, now define the following functions 
\[ \tilde{f}(x) \defeq \max_i \tilde{f}_i(x) \text{ where } \tilde{f}_i(x) \defeq f_i(x) + \left(1 - \frac{i}{T} \right)\gamma \defeq a_i^T x + \left(1 - \frac{i}{T} \right)\gamma\] 
\begin{equation}
  \label{eqn:hardfdef}
  \hardf(k, \gamma, \delta_T) \defeq S^k_{\delta_T, A_T} \;\tilde{f}
\end{equation}
i.e. smoothing $\tilde{f}$ with respect to $\delta_T, A_T$. The hard function we propose is the random function $\hardf$ with parameters set as $\gamma = \frac{1}{3\sqrt{T}}$ and $\delta_T = \frac{1}{20kT^{1.5}}$. 
  We restate facts which can be derived analogously to those derived in Section \ref{sec:construction} (c.f. Equations \eqref{eqn:errorbound},\eqref{eqn:funcminbound}).
  \begin{equation}
  \label{eqn:errorboundrand}
  \forall x \;\;|\hardf(x) - \tf(x)| \leq k\delta \;\;\text{ and }\;\; \min_{x \in B_d} \hardf(x) \leq \frac{-1}{\sqrt{T}} + \gamma + k\delta
  \end{equation}
The following key lemma will be the main component of the proof. 
\begin{lemma}
\label{lemma:randmainlemma}
  Let $\{x_1 \ldots x_T\}$ be the points queried by a randomized algorithm throughout its execution on the function $\hardf$. With probability at least $1 - \delta$ (over the randomess of the algorithm and the selection of $\hardf$) the following event $\mathcal{E}$ happens
  \[ \mathcal{E} = \left\{ \forall i \in [T] \;\; \forall j \geq i\;\;|a_j^Tx_i| \leq \smallq \right\} \]  
\end{lemma}
Using the above lemma we first demonstrate the proof of Theorem \ref{thm:mainthmrandomized}. We will assume the event $\mathcal{E}$ in Lemma \ref{lemma:randmainlemma} happens.
\\
\\

\noindent \textbf{Bounded Lipschitz Constants} Using Corollary \ref{cor:maincor}, the fact that $\hardf$ has Lipschitz constant bounded by 1 and that $\tilde{f}$ is invariant with respect to the $T$ dimensional subspace $A_T$, we get that the function $\hardf_T$ has higher order Lipschitz constants bounded above as
\[ \forall i \leq k \;\;\; L_{i+1} \leq \left(\frac{T}{\delta_T}\right)^i \leq \left(20kT^{2.5}\right)^i\] 

\noindent\textbf{Sub-optimality} : The event $\mathcal{E}$ in the lemma implies that $\tilde{f}_i(x_i) \geq -\smallq + (1 - \frac{i}{T})\gamma$ which implies that $\tilde{f}(x_i) \geq -\smallq + (1 - \frac{i}{T})\gamma$ and from Equation \eqref{eqn:errorboundrand} we get that 
\[ \forall i \in [T]\;\; \hardf(x_i) \geq -\smallq + \left(1 - \frac{i}{T}\right)\gamma - k\delta_T \]
Now using Equation \eqref{eqn:errorboundrand} we get that every $x_i$ is such that
\[ \forall i \in [T]\;\;\;\; \hardf(x_i) - \min_{x \in B} \hardf(x) \geq \frac{1}{\sqrt{T}} - \frac{i}{T}\gamma - 2k\delta_T - \smallq \geq \frac{1}{2\sqrt{T}}\]
The last inequality follows by the choice of parameters.
This finishes the proof of Theorem \ref{thm:mainthmrandomized}.
\end{proof}
\begin{proof}[Proof of Lemma \ref{lemma:randmainlemma}]
  We will use the following claims to prove the lemma. For any vector $x$, define the event $\mathcal{E}_i(x) = \big\{ \forall j \geq i\;\;|a_j^Tx| \leq \smallq \big\}$. The event we care about then is 
  \[ \mathcal{E} \defeq \bigcup_{i = 1 \rightarrow T} \mathcal{E}_i(x_i)\]

  \begin{claim}
  \label{claim:splineclaim}
      $\forall x$ if $\mathcal{E}_i(x)$ holds, then $[\hardf(x), \nabla \hardf(x) \ldots \nabla^k \hardf(x)]$ all depend only on $\{ a_1 \ldots a_i\}$. 
    \end{claim}
  \begin{claim}
  \label{claim:inductiveclaim}
      For any $i \in [T]$, if $\forall j < i$ if we have that $\mathcal{E}_j(x_j)$ holds then we have that with probability at least $1 - \frac{\delta}{T}$(over the choice of $a_i$ and the randomness of the algorithm) the event $\mathcal{E}_{i}(x_{i})$ happens.
    \end{claim}
    Claim \ref{claim:splineclaim} is a robust version of the argument presented in the proof of Theorem \ref{thm:mainthm}. Claim \ref{claim:inductiveclaim} is a byproduct of the fact that in high dimensions the corelation between a fixed vector and a random small basis is small. Claim \ref{claim:splineclaim} is used to prove the Claim \ref{claim:inductiveclaim}.

  Lemma \ref{lemma:randmainlemma} now follows via a simple inductive argument using Claim \ref{claim:inductiveclaim} which is as follows
      \[ \pr(\mathcal{E}) = \pr \left(\bigcup_{i = 1 \rightarrow T} \mathcal{E}_i(x_i)\right) = \prod_{i=1 \rightarrow T}\pr\left(\mathcal{E}_i(x_i) \big| \bigcup_{j < i} \mathcal{E}_j(x_j)\right) \geq \left(1 - \frac{\delta}{T}\right)^T \geq 1 - \delta\] 
\end{proof}
\begin{proof}[Proof of Claim \ref{claim:splineclaim}]
   As noted before the smoothing operator $S^k_{\delta}$ is such that at any point $x$ all the $k$ derivatives of $S^k_{\delta} f$ depend only on the value of $f$ in a ball of radius $k \delta$ around the point $x$. Therefore  it is sufficient to show that for the function $\tf$, for every $y$ such that $\|y - x\| \leq k\delta_T$ we have that $\tf(y)$ depends only on $\{a_1 \ldots a_{i}\}$. To ensure this, it is enough to ensure that for every such $y$ we have that $\argmin_{j \in [T]} \tf_j(y) \leq i$ which is what we prove next.

Lets first note the following facts. By the definition of $\mathcal{E}_i(x)$ we have that $\forall j > i, f_j(x) \leq \smallq$. This immediately implies that             
\begin{equation}
\max_{j > i} \;\; \tf_j(x_i) \leq \smallq + \left(1 - \frac{i+1}{T} \right) \gamma
\end{equation}
Now since we know each $\tf_i$ is $1$-Lipschitz \footnote{$\|a_i\| = 1$}, this also gives us 
\begin{equation}
\forall y \;\;\text{s.t.}\;\;\|y - x\| \leq k\delta_T \;\;\text{we have}\;\;\max_{j > i} \;\; \tf_j(y) \leq \smallq + \left(1 - \frac{i+1}{T} \right) \gamma + k\delta_T 
\end{equation}
By the event $\mathcal{E}_i(x)$ we also know that $\tf_i(x) \geq -\smallq + (1 - \frac{i}{T})\gamma$.
This implies as above 
\begin{equation}
\forall y \;\;\text{s.t.}\;\;\|y - x\| \leq k\delta_T \;\;\text{we have}\;\;\max_{j \leq i} \;\; \tf_j(y) \geq - \smallq + \left(1 - \frac{i}{T} \right) \gamma - k\delta_T 
\end{equation}
The above equations imply that as long as $2k\delta_T + \frac{1}{10T^{1.5}} < \frac{\gamma}{T}$ (which is true by the choice of parameters), we have that  
\begin{equation}
\forall y \;\;\text{s.t.}\;\;\|y - x\| \leq k\delta_T \;\; \argmin_{j \in [t]} \tf_j(y) \leq i 
\end{equation}
which is sufficient to prove Claim \ref{claim:splineclaim}.
 \end{proof} 
\begin{proof}[Proof of Claim \ref{claim:inductiveclaim}]
  Consider any $i \in [T]$. Given $\mathcal{E}_j(x_j)$ is true for all $j < i$, applying Claim \ref{claim:splineclaim} for all $j < i$, implies that all the information that the algorithm possesses is only a function of $\{a_1 \ldots a_{i-1}\}$ and the internal randomness of the algorithm. Further we can assume that the basis $\{a_1 \ldots a_T\}$ is chosen by the inductive process which picks $a_i$ uniformly randomly from the subspace $A_{i-1}^{\perp}$.

  Therefore we have that the choice of the remaining basis $\{a_i \ldots a_T\}$ is uniformly distributed in a subspace of dimension $d-i + 1$ completely independent of any vector $x_i$ the algorithm might play. Since we wish to bound the absolute value of the inner product we can assume $\|x_i\| = 1$ \footnote{Otherwise the absolute value of the inner product is only lower}.

  The lemma now reduces to the following quantity, consider a fixed unit vector $y$ in a $\reals^{d-i+1}$. Consider picking a $T - i + 1$ dimensional subspace of $\reals^{d - i + 1}$ given by the basis $y_1 \ldots y_{T-i+1}$ uniformly randomly. We wish to bound the probability that 
  \[Pr\left( \forall \;j\;\;|\innerprod{y}{y_j}| > \smallq\right) .\]
  The rest of the argument follows the argument by \cite{WoodworthSrebro2016}(Proof of Lemma 7). Note that 
  for $y_1$ this probability amounts to the surface area of a sphere above the caps of radius $\sqrt{1 - (\smallq)^2}$ relative to the surface area of the unit sphere. This surface area is smaller than the relative surface area of a sphere of radius $\sqrt{1 - (\smallq)^2}$. Formally this gives us 
  \[ Pr\left(|y_1^T y| \geq \frac{1}{20T^{1.5}}\right) \leq \sqrt{\left(1 - \left(\smallq\right)^2\right)^{d-i+1}}\]
  Applying the argument inductively we get that
  \[ \forall j \in [1, T-i+1] \;\;\; Pr\left(|y_j^T y| \geq \frac{1}{20T^{1.5}}\right) \leq \sqrt{\left(1 - \left(\smallq\right)^2\right)^{d-i-j+2}}\]
  Using the union bound we have that
  \begin{multline}
    Pr(\mathcal{E}_i(x_i)) \geq 1 - \bigcup_{j=1 \rightarrow T-i+1}\left( Pr\left(|y_j^T y| \geq \frac{1}{20T^{1.5}}\right)\right) \geq 1 - (T - i)\left(1 - \left(\smallq\right)^2\right)^{\frac{d-T+1}{2}} \\ \geq 1 - (T - i)e^{- (\smallq)^2\frac{d-T}{2}} \geq 1 - \frac{\delta}{T}
  \end{multline}
   \[ \] 
   The last line follows from the choice of $d = \Omega\left(T^3\log(\delta T^2)\right)$.
\end{proof}

\section{Acknowledgements}

The authors would like to acknowledge and thank Ohad Shamir for providing insightful comments on the first draft of this manuscript and Brian Bullins and Gopi Sivakanth for helpful suggestions. The first author is supported in part by Elad Hazan's NSF grant 1523815.

\bibliographystyle{alpha}
\bibliography{references}

\ENDDOC